\documentclass{article}
\usepackage{amssymb}
\usepackage{amsmath}
\usepackage{amsfonts}

\newtheorem{theorem}{Theorem}

\newtheorem{conjecture}[theorem]{Conjecture}

\input{tcilatex}

\begin{document}

\title{Cumulants of products of Normally distributed random variables}
\date{}
\author{Clarence KALITSI and Jan \ VRBIK \\
Department of Mathematics\\
Brock University, 500 GLenridge Ave.\\
St. Catharines, Ontario, Canada, L2S 3A1}
\maketitle

\begin{abstract}
To find moments of various estimators related to Autoregressive models of
Statistics, one first needs the cumulants of \emph{products} of two Normally
distributed random variables. The purpose of this article is to derive the
corresponding formulas, and extend them to products of three or more such
variables.
\end{abstract}

\section{Introduction}

The formulas presented in this article are crucial for finding moments (and,
subsequently, approximate distributions) of various parameter estimators
related to \textsf{AR(k)} models (see \cite{vrbik} and \cite{LIU}).

\subsection{Multivariate Normal distribution}

Assume that $X_{1},$ $X_{2},$... are \emph{centralized} (having zero mean),
Normally distributed random variables. Their moment-generating function is
given by%
\begin{equation}
\exp \left( \frac{\mathbf{X}^{\text{T}}\mathbb{V}~\mathbf{X}}{2}\right) 
\label{M}
\end{equation}%
where $\mathbb{V}$ is the corresponding variance/covariance matrix. Based on
(\ref{M}), we can easily find the expected value of a product of any number
of such random variables (which defines the corresponding \textsc{moment}),
getting \emph{zero} when this number is \emph{odd}, and%
\begin{equation}
\mathbb{E}\left( X_{1}X_{2}...X_{2k}\right) \equiv \mu
_{\{1,2,...2k\}}=\sum_{\mathbf{i}\epsilon \mathcal{P}%
_{2k}}V_{i_{1},i_{2}}V_{i_{3},i_{4}}...V_{i_{2k-1},i_{2k}}  \label{EV}
\end{equation}%
when this number is \emph{even}. The summation is over all possible $\binom{%
2k}{2,2,...2}/k!$ ways of \textsc{pairing} the $n$ indices, regardless of
the order (both within and between) of the resulting pairs. Thus, for
example, $\{1,2\},\{3,4\}$ and $\{4,3\},\{2,1\}$ are considered identical
and appear only once, but $\{1,2\},\{3,4\}$ and $\{1,3\},\{2,4\}$ are
distinct.

Notationally (and rather symbolically), this has been indicated by%
\begin{equation*}
\mathbf{i}\epsilon \mathcal{P}_{2k}
\end{equation*}%
where $\mathbf{i}$ represents the $\{i_{1},i_{2}\},$ $\{i_{3},i_{4}\},...$
indices and $\mathcal{P}_{2k}$ the set of all such selections.

\begin{description}
\item[Example] 
\begin{equation}
\mu _{\{1,2,3,4\}}=\mathbb{E}\left( X_{1}X_{2}X_{3}X_{4}\right)
=V_{1,2}V_{3,4}+V_{1,3}V_{2,4}+V_{1,4}V_{2,3}\ \ \ \ \ \ \ \ \square
\label{ex}
\end{equation}
\end{description}

Note that the resulting formulas are fully general; (\ref{ex}) can be read as%
\begin{equation*}
\mu
_{\{j_{1},j_{2},j_{3},j_{4}%
\}}=V_{j_{1},j_{2}}V_{j_{3},j_{4}}+V_{j_{1},j_{3}}V_{j_{2},j_{4}}+V_{j_{1},j_{4}}V_{j_{2},j_{3}}
\end{equation*}%
but using specific integers simplifies the notation; also, the formulas
allow any duplication of indices, e.g.%
\begin{equation*}
\mu _{\{1,1,2,2\}}=V_{1,1}V_{2,2}+2V_{1,2}^{2}
\end{equation*}

It is important to realize that, for random variables with zero mean (the
variables of this section), there is no difference between simple and
central moments; our $\mu $ thus stands for either.

\subsection{Joint cumulants}

Consider a collection of random variables (not necessarily centralized, nor
Normally distributed), say $Y_{1},$ $Y_{2},$..., and their joint
moment-generating function, defined by%
\begin{equation}
M(t_{1},t_{2},...)=\mathbb{E}\left[ \exp (t_{1}Y_{1}+t_{2}Y_{2}+...)\right]
\label{lnM}
\end{equation}%
The corresponding joint \textsc{cumulant} of $Y_{1},$ $Y_{2},...Y_{\ell }$
is defined by%
\begin{equation}
\kappa _{1,2,...\ell }\equiv \left. \frac{\partial ^{\ell }}{\partial
t_{1}\partial t_{2}...\partial t_{\ell }}\ln \left[ M(t_{1},t_{2},...)\right]
\right\vert _{t_{1}=t_{2}=...=t_{\ell }=0}  \label{CGF}
\end{equation}

It is well known (and easy to derive) that $\kappa _{i}=\mathbb{E}\left(
Y_{i}\right) ,$ and that all higher-order cumulants can be expressed in
terms of moments (of the $Y$ variables) thus:%
\begin{equation}
\kappa _{1,2,...\ell }=\sum_{\mathbf{j}_{1},\mathbf{j}_{2}...\mathbf{j}%
_{i}\epsilon \mathcal{A}(1,2,...\ell )}(-1)^{i-1}(i-1)!\mu _{\mathbf{j}%
_{1}}\mu _{\mathbf{j}_{2}}...\mu _{\mathbf{j}_{i}}  \label{cum}
\end{equation}%
where $\mathcal{A}(1,2,...\ell )$ is the collection of all \emph{partitions}
of the $1,2,...\ell $ indices. A \textsc{partition} is a division of a set
into an arbitrary number of non-empty and non-overlapping subsets (these are
denoted $\mathbf{j}_{1},\mathbf{j}_{2},\mathbf{...}$). For example, $%
\mathcal{A}(1,2,3)$ consists of $\mathbf{j}_{1}=\{1,2,3\},$ $\mathbf{j}_{1},%
\mathbf{j}_{2}=\{1\},\{2,3\},$ $\mathbf{j}_{1},\mathbf{j}_{2}=\{2\},\{1,3\},$
$\mathbf{j}_{1},\mathbf{j}_{2}=\{3\},\{1,2\},$ and $\mathbf{j}_{1},\mathbf{j}%
_{2},\mathbf{j}_{3}=\{1\},\{2\},\{3\}.$\newpage 

There are two points to make about (\ref{cum}):

\begin{itemize}
\item the formula is correct regardless of whether the moments are central
or simple (from now on, we denote \emph{central} moments by $\mu _{\mathbf{j}%
}$ and \emph{simple} moments by $\tilde{\mu}_{\mathbf{j}},$ respectively),

\item using central moments simplifies the RHS substantially - any partition
containing a \emph{single} index can be omitted (the corresponding $\mu _{%
\mathbf{j}}$ is zero).
\end{itemize}

We will spell out explicitly the first few of these formulas, first using 
\emph{simple} moments

\begin{eqnarray*}
\kappa _{1,2} &=&\tilde{\mu}_{\{1,2\}}-\tilde{\mu}_{\{1\}}\tilde{\mu}_{\{2\}}
\\
\kappa _{1,2,3} &=&\tilde{\mu}_{\{1,2,3\}}-\tilde{\mu}_{\{1,2\}}\tilde{\mu}%
_{\{3\}}-\tilde{\mu}_{\{1,3\}}\tilde{\mu}_{\{2\}}-\tilde{\mu}_{\{2,3\}}%
\tilde{\mu}_{\{1\}}+2\tilde{\mu}_{\{1\}}\tilde{\mu}_{\{2\}}\tilde{\mu}%
_{\{3\}} \\
&&\vdots
\end{eqnarray*}%
then, using \emph{central} moments (note the simplification):

\begin{eqnarray}
\kappa _{1,2} &=&\mu _{\{1,2\}}  \label{km} \\
\kappa _{1,2,3} &=&\mu _{\{1,2,3\}}  \notag \\
\kappa _{1,2,3,4} &=&\mu _{\{1,2,3,4\}}-\mu _{\{1,2\}}\mu _{\{3,4\}}-\mu
_{\{1,3\}}\mu _{\{2,4\}}-\mu _{\{1,4\}}\mu _{\{2,3\}}  \notag
\end{eqnarray}%
To continue, we consider only the special case of having all indices
identical (the corresponding general formulas get too lengthy):%
\begin{eqnarray*}
\kappa _{1,1,1,1,1} &=&\mu _{\{1,1,1,1,1\}}-10\mu _{\{1,1,1\}}\mu _{\{1,1\}}
\\
\kappa _{1,1,1,1,1,1} &=&\mu _{\{1,1,1,1,1,1\}}-15\mu _{\{1,1,1,1\}}\mu
_{\{1,1\}}-10\mu _{\{1,1,1\}}^{2}+30\mu _{\{1,1\}}^{3} \\
&&...
\end{eqnarray*}

Similarly to moments, cumulants are \emph{fully symmetric} in the
permutation of indices, e.g. $\kappa _{1,2,3,4}$ is the same as $\kappa
_{4,1,3,2}$ etc.

We now proceed to derive explicit formulas for these cumulants when some of
the $Y$ variables are centralized, Normally distributed (the $X$'s of the
previous section - we call them \emph{`}\textsc{singlets}\emph{'}), and the
others are \emph{products of two} such $X$'s (\textsc{doublets}). This poses
a bit of a notational challenge; we will use simple indices for singlets,
two indices in parentheses for doublets. For example, $\kappa _{1,2,(3,4)}$
indicates a third-order cumulant of three random variables, $X_{1},$ $X_{2}$
and $X_{3}X_{4}.$

\section{Cumulants involving singlets and/or doublets}

It is well known and easy to derive, by combining (\ref{M}) and (\ref{CGF}),
that all cumulants involving only \emph{singlets} are equal to zero, with
the exception of 
\begin{equation*}
\kappa _{1,2}=V_{1,2}
\end{equation*}

For \emph{doublets}, one can derive (by combining (\ref{cum}) and (\ref{EV}%
), and using a routine, `brute-force' computation - see the Appendix), that%
\begin{eqnarray*}
\kappa _{(1,2)} &=&V_{1,2} \\
\kappa _{(1,2),(3,4)} &=&V_{1,3}V_{2,4}+V_{1,4}V_{2,3} \\
\kappa _{(1,2),(3,4),(5,6)}
&=&V_{2,3}V_{4,5}V_{1,6}+V_{2,4}V_{3,5}V_{1,6}+V_{2,3}V_{4,6}V_{1,5}+V_{2,4}V_{3,6}V_{1,5}+
\\
&&V_{2,5}V_{3,6}V_{1,4}+V_{2,6}V_{3,5}V_{1,4}+V_{2,5}V_{4,6}V_{1,3}+V_{2,6}V_{4,5}V_{1,3}
\end{eqnarray*}%
and%
\begin{equation}
\kappa _{(1,2),...(2k-1,2k)}=\sum_{\mathbf{i}\epsilon \mathcal{\tilde{P}}%
_{2k}}V_{i_{1}i_{2}}V_{i_{3}i_{4}}...V_{i_{2k-1}i_{2k}}  \label{double}
\end{equation}%
in general, where $\mathcal{\tilde{P}}_{2k}$ indicates the set of all index
pairings (the old $\mathcal{P}_{2k}$) which \emph{do not} contain any of the
original $\{1,2\}$ or $\{3,4\}$ or ... $\{2k-1,2k\}$ pairs. There is a
simple scheme for building $\mathcal{\tilde{P}}_{2k}$:

\begin{itemize}
\item Start with the original pairs, e.g. $\{1,2\},\{3,4\},\{5,6\},$

\item keeping the first pair fixed, go over all $(k-1)!$ permutations of the
remaining pairs, e.g. $\{1,2\},\{5,6\},\{3,4\},$

\item keeping the first pair fixed, go over all $2^{k-1}$ interchanges of
indices in the remaining pairs, e.g. $\{1,2\},\{6,5\},\{3,4\},$ etc. (four
of them in this case),

\item shift each resulting arrangement by one index, e.g. $%
\{2,6\},\{5,3\},\{4,1\}.$
\end{itemize}

One can thus see that the number of terms on the RHS of (\ref{double}) is $%
(k-1)!\times 2^{k-1},$ which equals to $2,$ $8,$ $48$ when $k=2,$ $3$ and $4$
respectively (a fast-growing sequence).

\subsection{Mixed cases}

Let us now investigate cumulants with a mixture of singlets and doublets.
The rules for computing these prove to be quite simple:\bigskip 

A cumulant involving

\begin{itemize}
\item one singlet (regardless of the number of doublets) equals to zero,

\item more than two singlets (and any number of doublets) is also equal to
zero.

\item When a cumulant contains two singlets, it equals to the cumulant in
which the singlets are replaced by one doublet, i.e.%
\begin{equation*}
\kappa _{1,2,(3,4),...(2k-1,2k)}=\kappa _{(1,2),(3,4),...(2k-1,2k)}
\end{equation*}
\end{itemize}

\section{Beyond doublets}

The same approach enables us to develop formulas for cumulants which may
also involve triplets, quadruplets, etc. These are not needed when dealing
with parameter estimation related to \textsf{AR(k)} models, but may have a
potential application elsewhere.

Thus, for example, odd-order cumulants involving only triplets are all equal
to zero; for even-orders we get%
\begin{eqnarray*}
&&\kappa _{(1,2,3),(4,5,6)}\overset{}{=}\mu _{\{1,2,34,5,6\}} \\
&&\kappa _{(1,2,3),(4,5,6),(7,8,9),(10,11,12)}\overset{}{=}\mu
_{\{1,2,3,4,5,6,7,8,9,10,11,12\}} \\
&&-\mu _{\{1,2,3,4,5,6\}}\mu _{\{7,8,9,10,11,12\}}-\mu _{\{1,2,3,7,8,9\}}\mu
_{\{4,5,6,10,11,12\}} \\
&&-\mu _{\{1,2,3,10,11,12\}}\mu _{\{4,5,6,7,8,9\}} \\
&&\vdots 
\end{eqnarray*}%
following (\ref{km}). Note that the first of these cumulants turns out to be
a sum of $15$ terms of the 
\begin{equation*}
V_{i_{1},i_{2}}V_{i_{3},i_{4}}...V_{i_{2k-1},i_{2k}}
\end{equation*}%
type, while the second one has already $9720$ terms!

\begin{conjecture}
It seems that each cumulant with all indices distinct (regardless of how
they are grouped) can be always expressed in a form of the RHS of (\ref%
{double}), where $\mathcal{\tilde{P}}_{2k}$ is a specific subset of $%
\mathcal{P}_{2k}$ (and $2k$ is the total number of indices, which must be
even for all non-zero cumulants).
\end{conjecture}

Similarly, this is how one evaluates the first few cumulants involving only
quadruplets:%
\begin{eqnarray*}
&&\kappa _{(1,2,3,4)}\overset{}{=}\mu _{\{1,2,3,4\}} \\
&&\kappa _{(1,2,3,4),(5,6,7,8)}\overset{}{=}\mu _{\{1,2,3,4,5,6,7,8\}}-\mu
_{\{1,2,3,4\}}\mu _{\{5,6,7,8\}} \\
&&\kappa _{(1,2,3,4),(5,6,7,8),(9,10,11,12)}\overset{}{=}\mu
_{\{1,2,3,4,5,6,7,8,9,10,11,12\}} \\
&&-\mu _{\{1,2,3,4,5,6,7,8\}}\mu _{\{9,10,11,12\}}-\mu
_{\{1,2,3,4,9,10,11,12\}}\mu _{\{5,6,7,8\}} \\
&&-\mu _{\{5,6,7,8,9,10,11,12\}}\mu _{\{1,2,3,4\}}+2\mu _{\{1,2,3,4\}}\mu
_{\{5,6,7,8\}}\mu _{\{9,10,11,12\}}
\end{eqnarray*}%
resulting in $3,$ $96$ and $9504$ terms, respectively.

One can then deal with `mixed' cumulants in the same manner. There is of
course no limit as to their total number; we will give just one example:%
\begin{eqnarray*}
&&\kappa _{1,2,(3,4,5),(6,7,8),(9,10,11,12)}\overset{}{=}\mu
_{\{1,2,3,4,5,6,7,8,9,10,11,12\}}-\mu _{\{1,2\}}\mu
_{\{3,4,5,6,7,8,9,10,11,12\}} \\
&&-\mu _{\{1,3,4,5\}}\mu _{\{2,6,7,8,9,10,11,12\}}-\mu _{\{2,3,4,5\}}\mu
_{\{1,6,7,8,9,10,11,12\}}-\mu _{\{1,6,7,8\}}\mu _{\{2,3,4,5,9,10,11,12\}} \\
&&-\mu _{\{2,6,7,8\}}\mu _{\{1,3,4,5,9,10,11,12\}}-\mu
_{\{1,2,3,4,5,6,7,8\}}\mu _{\{9,10,11,12\}}-\mu _{\{3,4,5,6,7,8\}}\mu
_{\{1,2,9,10,11,12\}} \\
&&+2\mu _{\{1,2\}}\mu _{\{3,4,5,6,7,8\}}\mu _{\{9,10,11,12\}}+2\mu
_{\{1,3,4,5\}}\mu _{\{2,6,7,8\}}\mu _{\{9,10,11,12\}} \\
&&+2\mu _{\{2,3,4,5\}}\mu _{\{1,6,7,8\}}\mu _{\{9,10,11,12\}}
\end{eqnarray*}%
resulting in a sum of $7848$ terms.

\section{Appendix}

For readers familiar with Mathematica, we now supply a few Mathematica
functions to facilitate the computation of cumulants involving centralized,
Normally distributed random variables and their products.

\subsection{Finding moments}

Computing the mean value (MV) of a product of centralized, Normally
distributed random variables can be achieved by:\bigskip

\noindent MV[a\_\_] := Module[\{a0 = Sort[Flatten[\{a\}]], a1, a2, a3\},

a1 = Length[a0]; a2 = Table[\{a0[[1]], a0[[i]]\}, \{i, 2, a1\}];

a3 = Table[Delete[a0, \{\{1\}, \{i\}\}], \{i, 2, a1\}];

Which[Mod[a1, 2] == 1, 0, a1 == 2, Apply[V, a0], True,

Apply[Plus, Table[Apply[V, a2[[i]]] MV[a3[[i]]], \{i, a1 - 1\}]]]//
Expand]\bigskip

To use it, we type\bigskip

\noindent MV[2, 5, 2, 5, 2, 8]/.V[i\_, j\_]-\TEXTsymbol{>}V$_{\text{i,j}}$%
\bigskip

\noindent where the integers in square brackets represent the corresponding
indices, listed in any order and allowing for any amount of duplication.

This returns the following answer:%
\begin{equation*}
3V_{2,2}V_{5,5}V_{2,8}+6V_{2,5}^{2}V_{2,8}+6V_{2,2}V_{2,5}V_{5,8}
\end{equation*}

\subsection{Computing cumulants}

This time our goal is a bit more ambitious: we want to build a Mathematica
program for computing cumulants of any number of singlets, doublets,
triplets, etc. of centralized, Normally distributed random variables.

We start by constructing a few auxiliary Mathematica functions:\bigskip

\noindent aux1[b\_] := Module[\{b1 = Union[b], b2, b3\},

b2 = Map[Position[b, \#][[1, 1]] \&, b1];

b3 = Length[b2]; Table[\{b[[b2[[i]]]], Delete[b, b2[[i]]]\}, \{i,
b3\}]]\bigskip

\noindent aux2[a\_, B\_] := Module[\{b1 = Map[Join[\{a[[1]]\}, \#] \&,

KSubsets[Drop[a, 1], B[[1]] - 1]]\},

Flatten[Map[aux3[a, \#, B[[2]]] \&, b1], 1]]\bigskip

\noindent aux3[a\_, b\_, B\_] := Map[Join[\{b\}, \#] \&, aux4[Complement[a,
b], B]]\bigskip

\noindent aux4[a\_, b\_] := Module[\{B = aux1[b]\}, If[Length[b] == 1,
\{\{a\}\},

Apply[Join, Table[aux2[a, B[[i]]], \{i, Length[B]\}]]]]\bigskip

The program to compute a cumulant of any number of products then looks like
this:\bigskip 

\noindent K[a\_\_] := Module[\{a0 = \{a\}, L, q, w\}, L = Length[a0];

q = Map[Reverse, Partitions[L]];

w = Map[aux4[Range[L], \#] \&, q] /. Table[i -\TEXTsymbol{>} a0[[i]], \{i,
L\}];

Apply[Plus, Apply[Times, Apply[MV, w, \{3\}], \{2\}], \{1\}] .

Map[(-1)$^{\text{Length[\#]-1}}$(Length[\#] - 1)! \&, q] // Expand]\bigskip

To find the fifth-order cumulant of the random variables $X_{3},$ $%
X_{1}X_{3},$ $X_{1}X_{3},$ $X_{1}X_{2}X_{3}$ and $X_{1}X_{2}X_{3}^{2},$ one
has to type (to simplify the answer, we have assumed that the $X$'s are 
\emph{standardized}, i.e. have a mean of zero \emph{and} the variance equal
to $1$):\bigskip

K[3,\{1,3\},\{1,3\},\{1,2,3\},\{1,2,3,3\}]/.V[i\_,i\_]-\TEXTsymbol{>}%
1/.V[i\_,j\_]-\TEXTsymbol{>}C$_{\text{i,j}}\bigskip $

\noindent resulting in%
\begin{eqnarray*}
&&42+158C_{1,2}^{2}+438C_{1,3}^{2}+240C_{2,3}^{2}+1784C_{1,2}C_{1,3}C_{2,3}+1960C_{1,2}C_{1,3}^{2}C_{2,3}
\\
&&+802C_{1,2}^{2}C_{1,3}^{2}+1616C_{1,3}^{2}C_{2,3}^{2}+240C_{1,3}^{4}+400C_{1,3}^{4}C_{2,3}^{2}
\end{eqnarray*}%
where $C_{i,j}$ is the correlation coefficient between $X_{i}$ and $X_{j}.$
Note that\bigskip

K[\{3,1\},\{3,2,1\},3,\{2,3,3,1\},\{1,3\}]/.V[i\_,i\_]-\TEXTsymbol{>}%
1/.V[i\_,j\_]-\TEXTsymbol{>}C$_{\text{i,j}}\bigskip $

\noindent would have returned the same answer.

\end{document}